\documentclass[a4paper, 12pt]{article}
\usepackage{amssymb, amsmath, amsfonts}
\usepackage{amsthm,amscd,amsfonts,latexsym,color,epsfig}
\begin{document}

\sloppy 

\begin{center}
{\Large \bf Inverse source problems for time-fractional mixed parabolic-hyperbolic type equations}

\medskip

{\bf Pengbin Feng}

{\it Academy of Mathematics and Systems Science, Chinese Academy of Sciences, Zhongguancun East Road 55, 100190 Beijing, P.R.China}\\
{fengpengbin11@mails.ucas.ac.cn}

{\bf E.T.Karimov}

{\it Institute of Mathematics, National University of Uzbekistan, Durmon yuli str.,29, 100125 Tashkent, Uzbekistan}\\
{erkinjon@gmail.com}

\end{center}

{\bf Abstract}

In the present paper we consider an inverse source problem for time-fractional mixed parabolic-hyperbolic equation with the Caputo derivative. In case, when hyperbolic part of the considered mixed type equation is wave equation, the uniqueness of source and solution are strongly influenced by initial time and generally is ill-posed. However, when the hyperbolic part is time fractional, the problem is well posed if end time is large. Our method relies on the orthonormal system of eigenfunctions of the operator with respect to space variable. We proved the uniqueness and stability of certain weak solutions for considered problems.

\section{Introduction}

Theory of boundary problems for fractional order differential equations is one of the rapidly developing branches of the Fractional Calculus. Since many mathematical models of real-life processes are directly connected with the investigations on aforementioned theory, it becomes very popular among the specialists on differential equations.

Omitting many papers on direct boundary problems for PDEs involving fractional differential operators, we note some works [1-3], where time-fractional parabolic-hyperbolic type equations were investigated. Precisely, main boundary problems for mixed parabolic-hyperbolic equations with the Riemann-Liouville fractional differential operator in parabolic part, were objects of investigations. In aforementioned works "hyperbolic part" of the considered mixed domain is characteristic triangle. Therefore, used methods are different than method we would like to handle in the present work. Hence, obtained results vary as well.

We would like to note growing interest of specialists to the inverse problems for fractional order PDEs. Especially, inverse problems related to the finding space-dependent source using the additional boundary measurement, become interesting. For instance, in [4] such problem with non-local conditions with respect to space variable for time-fractional diffusion equation was investigated. Due to non-local conditions, authors express seeking functions by special bi-orthoghonal series. 

Results of the present paper is closely related to the work [5] and contain its results as particular case (see Remark 3.6). Obtained results could be useful for investigations on parabolic-hyperbolic equations, using numerical methods, for example, see [6].

Regarding the inverse problems for time-fractional diffusion-wave equations, considering uniformly elliptic operator in the space variables, we refer readers to the works [7-9].

We would like to note that solvability of boundary problems for mixed type equations directly depends from so-called "gluing conditions". On the line of type changing we need to glue value of seeking function and value of its derivative in order to get solution in a whole domain. There exist many types of gluing conditions such as continuous, discontinuous, integral form and etc. For instance, in the works [2,3] gluing conditions of integral form were in use, but in the work [5] authors consider boundary problems with continuous gluing conditions, i.e values of seeking function and its derivative from the both parabolic and hyperbolic parts of mixed domain are equal on the line of type changing. Depending on which gluing conditions are used, solvability conditions to given data vary. Some physical meaning of gluing conditions for parabolic-hyperbolic type equations one can find in the monograph [10].

In the present work, due to time-fractional parabolic equation, we used special gluing condition, which depends from the fractional order $\alpha$. In particular integer case, i.e. $\alpha=1$, we will get continuous gluing condition.

The rest of this paper is organized as follows. In section 2, we give some preliminaries regarding to the definition of the Caputo fractional differential operator, general solution of fractional differential equation with the Caputo derivative and some properties of the Mittag-Leffler functions. In sections 3.1 and 3.2 we formulate problem for mixed type equation containing fractional diffusion and wave equation. Section 3.3 is devoted to the proof of the uniqueness of the solution. In section 3.4 we give existence and stability results. Finally, in section 4 we formulate another problem for purely time-fractional mixed parabolic-hyperbolic equation. We state the result on unique weak solvability of the problem.
In Appendix section one can find some verifications of statements.

\section{Preliminaries}
\subsection{Definition of the Caputo fractional differential operator and some properties of two parameter Mittag-Leffler function}
Below we give a definition of the Caputo fractional differential operator (see [11, p.14])
The expression
$$
{}_{C}D_{st}^{\alpha }\varphi \left( t \right)=sign^n(t-s) D_{st}^{\alpha -n}{{\varphi }^{\left( n \right)}}\left( t \right),\,n-1<\alpha \le n,\,\,\,n\in {\bf N}
$$
we call as the Caputo fractional differential operator of the order $\alpha $. For a function $\varphi \left( y \right)$, the Riemann-Liouville integral-differential operator of the order $\alpha $ with initial point $s\in$ {\bf R}, can be defined as follows:
$$
D_{st}^{\alpha }\varphi \left( t \right)\equiv \left\{ \begin{array}{l}
\frac{sign(t-s)}{\Gamma \left( -\alpha  \right)}\int\limits_{s}^{t}{\frac{\varphi \left( z \right)dz}{{{\left| t-z \right|}^{\alpha +1}}},\,\,\alpha <0,} \\
\varphi \left( t \right),\,\,\,\,\,\,\,\,\,\,\,\,\,\,\,\,\,\,\,\,\,\,\,\,\,\,\,\,\,\,\,\alpha =0, \\
sign^n(t-s)\frac{{{d}^{n}}}{d{{t}^{n}}}D_{st}^{\alpha -n}\varphi \left( t \right),\,\,\,\,\,\,\,\,\,\,\,n-1<\alpha \le n,\,n\in {\bf N}. \\
\end{array} \right.
$$

The Mittag-Leffler function of two parameter is defined as [12, p.17]
$$
E_{\alpha,\beta}(z)=\sum_{n=0}^\infty\frac{z^n}{\Gamma(\alpha n+\beta)},\,\,\,\alpha>0,\,\beta>0.\eqno (2.1)
$$
We use the formula for the derivative of this function [12, p.21]:
$$
D_{0t}^\gamma\left(t^{\alpha k+\beta-1}E_{\alpha,\beta}^{(k)}\left(\lambda t^\alpha\right)\right)=t^{\alpha k+\beta-\gamma-1} E_{\alpha,\beta-\gamma}^{(k)}\left(\lambda t^\alpha\right),\eqno (2.2)
$$
where $\gamma$ is any arbitrary real number, $E_{\alpha,\beta}^{(k)}(y)=\frac{d^k}{dy^k}E_{\alpha,\beta}(y)$.

We, as well, will use the following property of this function [13, p.45]:
$$
E_{\alpha,\mu}(z)=\frac{1}{\Gamma(\mu)}+zE_{\alpha,\mu+\alpha}(z).\eqno (2.3)
$$
We need two asymptotic expansions of the Mittag-Leffler function, given below as theorems [12, p.35].

{\bf Theorem 2.1-1.} If $\alpha<2$, $\beta$ is arbitrary real number, $\mu$ is such that $\pi\alpha/2<\mu<\min \{\pi,\pi\alpha\}$ and $C_1$ and $C_2$ are real constants, then
$$
\left|E_{\alpha,\beta}(z)\right|\leq C_1 (1+|z|)^{(1-\beta)/\alpha}\exp \left(\Re\left(z^{1/\alpha}\right)\right)+\frac{C_2}{1+|z|},
\,\left(|\arg(z)|\leq \mu\right),\, |z|\geq 0.
$$

{\bf Theorem 2.1-2.} If $\alpha<2$, $\beta$ is arbitrary real number,
$\mu$ is such that $\pi\alpha/2<\mu<\min \{\pi,\pi\alpha\}$ and $C$ is real constant, then
$$
\left|E_{\alpha,\beta}(z)\right|\leq \frac{C}{1+|z|},\,\,\left(\mu\leq|\arg(z)|\leq \pi\right), \,|z|\geq 0.
$$

{\bf Theorem 2.2.} If $0<\alpha<2$, $\beta$ is arbitrary real number,
$\mu$ is such that $\pi\alpha/2<\mu<\min \{\pi,\pi\alpha\}$ then for any arbitrary integer $p\ge1$, the following expansion holds:
$$
E_{\alpha,\beta}(z)=\sum_{k=1}^{p} \frac{z^{-k}}{\Gamma(\beta-\alpha k)}+O(\left|z\right|^{-1-p}), \,\,\, \left(\mu\leq|\arg(z)|\leq \pi\right),\,\, \left|z\right|\to \infty.
$$

\subsection{Solutions of fractional order differential equations with the Caputo time fractional derivative.}

For the reader's convenience we rewrite the following fact taken from the monograph [11, p.17]:

The following Cauchy problem:
$$
\left(_CD_{sx}^\alpha y\right)(x)-\lambda y(x)=g(x), \,\,(s\in {\bf R},\,n-1<\alpha\le n,\,n\in {\bf N}),
$$
$$
y^{(n-k)}(0)=b_k,\,\,(k=1,2,..., n)
$$
has the unique solution represented as
$$
y(x)=\sum_{k=1}^{n}b_k |x-s|^{n-k} E_{\alpha, n-k+1}\left(\lambda |x-s|^\alpha\right)+sign(s-x)\int\limits_s^x |x-t|^{\alpha-1}E_{\alpha,\alpha}\left[\lambda|x-t|^\alpha\right]g(t)dt.\eqno (2.4)
$$
Afore-mentioned solution can be found as well in the work [14].

\section{Case, when hyperbolic part of the mixed type equation is integer order wave equation}

\subsection{Formulation of a problem}
In a domain $\Omega=\left\{(x,t):\,a<x<b,\,-p<t<q\right\} $ we consider an equation
$$
f(x)=\left\{
\begin{array}{l}
_CD_{at}^\alpha u -\mathcal{L}u,\,t>0\hfill\\
u_{tt} -\mathcal{L}u,\,t<0,\hfill\\
\end{array}
\right.\eqno (3.1)
$$
together with conditions
$$
u(a,t)=0,\,\,u(b,t)=0,\,\,-p\leq t \leq q, \eqno (3.2)
$$
$$
u(x,-p)=\psi(x),\,\,u(x,q)=\varphi(x),\,\,a\le x\le b, \eqno (3.3)
$$
$$
u(x,+0)=u(x,-0),\,a\le x\le b,\,\,\lim_{t\rightarrow +0} {_CD_{at}^\alpha u(x,t)}=u_t(x,-0),\,\,a<x<b,\eqno (3.4)
$$
where $\mathcal{L}$ is operator in $L^2(a,b)$ defined as follows
$$
\left\{
\begin{array}{l}
\displaystyle{\mathcal{L}v(x)=\frac{d}{d x}\left(r(x)\frac{d v}{d x}\right)-e(x)v,\,a<x<b,}\hfill\\
\mathcal{D}\left(\mathcal{L}\right)=\left\{v\in H^2(a,b);\,v(a)=v(b)=0\right\},\hfill\\
\end{array}
\right.
\eqno (3.5)
$$
$0<\alpha\le 1$, $p,q>0$, $a,b\in {\bf R}$, $\psi(x)$, $\varphi(x)$, $r(x)$, $e(x)$ are given functions such that $r(x)\in C^2[a,b],\,r(x)>0$, $e(x)\in C[a,b]$.

\subsection{Reformulation of a problem and formal construction of the solution}

Let us consider equation $\mathcal{L}v=-\mu v$. Using the Liouville transformation we have
$$
\bar{v}=l(x(z))v(x(z)),\,\,z\in[0,1],
$$
where
$$
l(x)=\left(r(x)\right)^{1/4},\,z(x)=\frac{1}{K}\int\limits_a^x\sqrt{\frac{1}{r(s)}}ds,\,K=\int\limits_a^b\sqrt{\frac{1}{r(s)}}ds,\,s\in [a,b],
$$
and further we get
$$
\mathcal{L}^*\bar{v}=-\lambda \bar{v},
$$
in which
$$
\left\{
\begin{array}{l}
\displaystyle{\mathcal{L}^*\bar{v}=\frac{d^2}{dz^2}\bar{v}(z)-g(z)\bar{v}(z),\,0<x<1,}\hfill\\
\displaystyle{\mathcal{D}^*\left(\mathcal{L}^*\right)=\left\{\bar{v}\in H^2(0,1);\,\bar{v}(0)=\bar{v}(1)=0\right\},}\hfill\\
\end{array}
\right.
\eqno (3.6)
$$
$$
\lambda=K^2\mu,\,g(z)=K^2\left[e(x(z))+l(x(z))\left(\frac{d(x(z))l'(x(z))}{l^2(x(z))}\right)'\right].
$$
Here $z(x)$ is monotone about $x$, and $x(z)$ is inverse function of $z(x)$. Note that $f(x),\,\varphi(x),\,\psi(x)$ as well were transformed to $\bar{f}(z),\, \bar{\varphi}(z),\,\bar{\psi}(z)$ respectively. Due to the properties of Liouville transformation, the uniqueness of $f(x)$ corresponds to the uniqueness of the $\bar{f}(z)$.

\smallskip

Further, for convenience we denote $\bar{v}(z),\, \bar{f}(z),\,  \bar{\varphi}(z),\,\bar{\psi}(z)$ as $u(x),\,f(x),\,\varphi(x),\,\psi(x)$ respectively.

According to this transformation, the domain $\Omega$ transferred to the\\ $\Omega^*=\left\{(x,t):\,0<x<1,\,-p<t<q\right\}$ and the problem (3.1)-(3.4) equivalently reduced to the following problem:
$$
f(x)=
\left\{
\begin{array}{l}
_CD_{0t}^\alpha u-\mathcal{L}^*u,\,t>0\hfill\\
u_{tt}-\mathcal{L}^*u,\,\,\,\,\,\,\,\,\,\,t<0,\hfill\\
\end{array}
\right.\eqno (3.7)
$$
$$
\left\{
\begin{array}{l}
u(0,t)=u(1,t)=0,\,-p\le t\le q,\hfill\\
u(x,-p)=\psi(x),\,u(x,q)=\varphi(x),\,0\le x\le 1,\hfill\\
\displaystyle{u(x,+0)=u(x,-0),\,0\le x\le 1,\,\,\lim_{t\rightarrow +0} {_CD_{0t}^\alpha u(x,t)}=u_t(x,-0),\,\,0<x<1.}\hfill\\
\end{array}
\right.\eqno (3.8)
$$

First we give the definition of the weak solution.

{\bf Definition 3.1.} We call $u(x,t)$ a \textit{weak solution} to (3.7)-(3.8), if $f(x)\in L^2(0,1)$ and $u(\cdot, t)\in\mathcal{D}^*\left(\mathcal{L}^*\right)\equiv H^2(0,1)\cap H_0^1(0,1)$ for $t\in [-p, q]$
and
$$
\begin{array}{l}
\displaystyle{u\in C\left([-p, q];\mathcal{D}^*\left(\mathcal{L}^*\right)\right),\,
\frac{\partial u}{\partial t}\in C\left([-p,0]\cup (0,q];L^2(0,1)\right),}\hfill\\
\displaystyle{{}_CD_{0t}^\alpha u\in C\left([0,q];L^2(0,1)\right),
\frac{\partial^2 u}{\partial t^2} \in C\left((-p,0);L^2(0,1)\right),}\hfill\\
\end{array}
\eqno (3.9)
$$
$$
\lim_{t\rightarrow -p}\left\|u(\cdot,t)-\psi\right\|_{\mathcal{D}^*\left(\mathcal{L}^*\right)}=0,\eqno (3.10)
$$
$$
\lim_{t\rightarrow q}\left\|u(\cdot,t)-\varphi\right\|_{\mathcal{D}^*\left(\mathcal{L}^*\right)}=0,\eqno (3.11)
$$
$$
\lim_{|t|\rightarrow 0}\left\| {_C}D_{0t}^\alpha u(\cdot,t)-u_t(\cdot,t)\right\|_{L^2}=0, \eqno (3.12)
$$
$$
\left({}_CD_{0t}^\alpha u,\eta\right)+\left(\mathcal{L}^*u,\eta\right)=\left(f,\eta\right),\,\, t\in [0,q],\,\forall \eta\in \mathcal{D}^*\left(\mathcal{L}^*\right),\eqno (3.13)
$$
$$
\left(\frac{\partial^2 u}{\partial t^2},\eta\right)+\left(\mathcal{L}^*u,\eta\right)=\left(f,\eta\right),\,\, t\in (-p,0),\,\forall \eta\in \mathcal{D}^*\left(\mathcal{L}^*\right). \eqno (3.14)
$$

{\bf Problem 1.} To determine uniquely a pare $\{u(x,t),f(x)\}$ in the domain $\Omega^*$, satisfying (3.7), (3.9)-(3.14).

Eigenvalues of the symmetric operator $\mathcal{L}^*$ are defined as $\lambda_k$ and corresponding complete systems of eigenfunctions as $\omega_k$. Moreover, since $g\in L^{\infty}[0,1]$ and is real-valued,  $\lambda_k\,(k=1,2,...)$ of the operator $\mathcal{L}^*$ are real-valued, simple and
$$
\left\{\lambda_1<\lambda_2<...<\lambda_l\le 0<\lambda_{l+1}...<\lambda_k<...\rightarrow +\infty\right\}.
$$ for simplicity, we assume $g\ge0$ in its domain, thus we have $l=0$ and all the eigenvalues are positive .
Asymptotic behavior is [15, p.135]
$$
\lambda_k=k ^2\pi^2+\int\limits_0^1g(x)dx-\int\limits_0^1g(x)\cos(2k\pi x)dx+O\left(\frac{1}{k}\right)\eqno (3.15)
$$
as $k\rightarrow \infty$ uniformly for $g\in L^\infty(0,1)\subset L^2(0,1)$.
Second integral of (3.15) is the $k$th Fourier coefficient of $g$ with respect to $\left\{\cos(2k\pi x):\,k=0,1,...\right\}$. Since $g(x)\in L^2(0,1)$, then this integral tends to $0$ at $k\rightarrow+\infty $.
Furthermore, the corresponding eigenfunctions $\omega_k$,
normalized to $\left\|\omega_k\right\|_{L^2}=1$ and have the following asymptotic behavior:
$$
\omega_k(x)=\sqrt{2}\sin(k\pi x)+O(1/n),\,\,\omega_k'(x)=\sqrt{2}k\pi\cos(k\pi x)+O(1)
$$
as $k\rightarrow \infty$ uniformly for $x\in[0,1]$ and $g\in L^\infty(0,1)\subset L^2(0,1)$.

Solution of the problem 1 we search as follows
$$
u(x,t)=\sum_{k=1}^{\infty}V_k(t)\omega_k(x),\,\,t>0, \eqno (3.16)
$$
$$
u(x,t)=\sum_{k=1}^{\infty}W_k(t)\omega_k(x),\,\,t<0 \eqno (3.17)
$$
and right-hand side as
$$
f(x)=\sum_{k=1}^{\infty}f_k\omega_k(x),\eqno (3.18)
$$
where
$$
\varphi_k=(\varphi,\omega_k)\quad \psi_k=(\psi,\omega_k), \eqno (3.19)
$$
$(\cdot,\cdot)$ is a scalar product in $L^2(0,1)$.

By standard calculation, it easily reduce that:

$$
V_k(t)=V_k(0)E_{\alpha,1}\left(-\lambda_kt^\alpha\right)+f_k t^\alpha
E_{\alpha,\alpha+1}\left(-\lambda_kt^\alpha\right),\eqno (3.20)
$$
$$
W_k(t)=A_k\sin \sqrt{\lambda_k}t +B_k\cos \sqrt{\lambda_k} t +\frac{f_k}{\lambda_k}.\eqno(3.21)
$$
Besides by the conditions (3.8) we obtain,
$$
V_k(0)=B_k +\frac{f_k}{\lambda_k},\,\,
\psi_k=A_k\sin \left(-\sqrt{\lambda_k}p\right)+B_k\cos \left(-\sqrt{\lambda_k}p\right) +\frac{f_k}{\lambda_k},\eqno\\
$$
$$
\varphi_k=V_k(0)E_{\alpha,1}\left(-\lambda_kq^\alpha\right)+f_kq^\alpha E_{\alpha,\alpha+1}\left(-\lambda_kq^\alpha\right),\,\,
f_k-\lambda_k V_k(0)=\sqrt{\lambda_k}A_k.
$$
By all these conditions above, we actually get:
$$
A_k=-\sqrt{\lambda_k}\frac{\varphi_k-\psi_k}{\Delta_k},\,\,
B_k=\frac{\varphi_k-\psi_k}{\Delta_k},
$$
$$
V_k(0)=\frac{\left(\varphi_k-\psi_k\right)\left(1-E_{\alpha,1}\left(-\lambda_kq^\alpha\right)\right)}{\Delta_k}+\varphi_k,
$$
$$
f_k=-\frac{\lambda_k\left(\varphi_k-\psi_k\right)E_{\alpha,1}\left(-\lambda_kq^\alpha\right)}{\Delta_k}+\lambda_k\varphi_k,
$$

where we denote
$$\Delta_k=E_{\alpha,1}\left(-\lambda_k q^\alpha\right)-\sqrt{\lambda_k} \sin \sqrt{\lambda_k} p-\cos \sqrt{\lambda_k}p.\eqno(3.22)
$$

Then considering (3.16)-(3.18), (3.20), (3.21) the formal solution of the problem we represent as:
$$
u(x,t)=\sum_{k=1}^{\infty}\left\{\frac{\varphi_k-\psi_k}{\Delta_k}\left[E_{\alpha,1}\left(-\lambda_kt^\alpha\right)-E_{\alpha,1}\left(-\lambda_kq^\alpha\right)\right]+\varphi_k\right\}\omega_k(x),\,t\in [0,q],
\eqno (3.23)
$$
$$
u(x,t)=\sum_{k=1}^{\infty}\left\{\frac{\varphi_k-\psi_k}{\Delta_k}\left[-\sqrt{\lambda_k}\sin\sqrt{\lambda_k}t+\cos\sqrt{\lambda_k}t-E_{\alpha,1}\left(-\lambda_kq^\alpha\right)\right]+\varphi_k\right\}\omega_k(x), \eqno (3.24)
$$
$$
t\in [-p,0],
$$
and
$$
f(x)=\sum_{k=1}^{\infty}\left\{\frac{-\lambda_k\left(\varphi_k-\psi_k\right)}{\Delta_k}E_{\alpha,1}\left(-\lambda_kq^\alpha\right)+\lambda_k\varphi_k\right\}\omega_k(x),
\eqno (3.25)
$$
where $\varphi_k$, $\psi_k$ are defined by (3.19) and $\Delta_k$ by (3.22).
\subsection{Conditional uniqueness of the solution and ill-posedness.}

{\bf Theorem 3.1.} If $\Delta_k\ne 0$, then formal solution of the problem 1 is unique.

{\it Proof.} 

Let $u(x,t)$ and $f(x)$ be a solution to problem 1 with $\psi(x)=\varphi(x)=0$, denote
$$
u_k(t)=\left(u(x,t),\omega_k(x)\right).
$$

Applying fractional operator ${}_CD_{0t}^\alpha$ to both sides of the equality above, at $t\in [0,q]$, considering 
$$
\left({}_CD_{0t}^\alpha u(\cdot,t),\omega_k(x)\right)={}_CD_{0t}^\alpha \left(u(\cdot,t),\omega_k(x)\right),\,\,0\le t\le q,
$$
(See the proof of this statement in Appendix section or refer to the work [9, Lemma A.1].)
and taking into account the boundary condition (3.11) together with $\varphi(x)=0$, we have
$$
_CD_{0t}^\alpha u_k(t)+\lambda_ku_k(t)=f_k,\, t\in [0,q].
$$
Similarly, considering (3.10) at $\psi(x)=0$ for $t\in [-p,0]$, we obtain
$$
u_k''(t)+\lambda_ku_k(t)=f_k, t\in [-p,0],
$$
here $u_k(-p)=u_k(q)=0$ and $f_k=\left(f(x),\omega_k(x)\right)$.

Solving them, we deduce
$$
u_k(t)=C_kE_{\alpha,1}\left(-\lambda_k t^\alpha\right)+f_k t^\alpha E_{\alpha,\alpha+1}\left(-\lambda_k t^\alpha\right),\ \ \  t\in [0,q] ,
$$
$$
u_k(t)=C_1\sin \sqrt{\lambda_k} t +C_2\cos \sqrt{\lambda_k} t +\frac{f_k}{\lambda_k},\ \ \ \    t\in [-p,0].
$$
Considering gluing condition in (3.9) and (3.12), we deduce
$$
\left\{
\begin{array}{l}
C_2=C_k-\frac{f_k}{\lambda_k}, \\
C_1=-\sqrt{\lambda_k} C_k +\frac{f_k}{\sqrt{\lambda_k}},\\
\end{array}
\right.
$$
Getting all these, we have 
$$
\left\{
\begin{array}{l}
C_k E_{\alpha,1}\left(-\lambda_kq^\alpha\right)+f_k q^\alpha E_{\alpha,\alpha+1}\left(-\lambda_kq^\alpha\right)=0, \\
C_k \lambda_k\left(\cos \sqrt{\lambda_k}p+\sqrt{\lambda_k} \sin \sqrt{\lambda_k}p\right)+f_k\left(1-\sqrt{\lambda_k}\sin\sqrt{\lambda_k}p-\cos\sqrt{\lambda_k}p\right)=0,\\
\end{array}
\right.
$$
the determination of this equation is exactly (3.22).

Since we suppose that $\Delta_k\ne 0$, then we have $C_k=f_k=0$ for any $k$ ,which shows $u_k\equiv 0$, because of the completeness of the base $\{\omega_k(x),\,k\in {\bf N}\}$, we can state that $u\equiv0$ and $f\equiv0$.\\
Theorem 3.1. is proved.

{\bf Remark 3.1.} For infinite many $p\in R^+$, it easily shows $\Delta_k=0$, then the problem 1 is ill-posed as the following example shows.\\
{\bf Example:} $\Delta_k=0$ if and only if

$$
p=\left\{
\begin{array}{l}
\displaystyle{\frac{1}{\sqrt{\lambda_k}}\left[\arcsin \frac{E_{\alpha,1}\left(-\lambda_kq^\alpha\right)}{\sqrt{\lambda_k+1}}-\gamma_k\right]+\frac{2n\pi}{\sqrt{\lambda_k}},} \\
\displaystyle{-\frac{1}{\sqrt{\lambda_k}}\left[\arcsin \frac{E_{\alpha,1}\left(-\lambda_kq^\alpha\right)}{\sqrt{\lambda_k+1}}+\gamma_k\right]+\frac{(2n+1)\pi}{\sqrt{\lambda_k}},}\\
\end{array}
\right.
$$
where  $\displaystyle{\gamma_k=\arcsin \frac{1}{\sqrt{\lambda_k+1}}}$, $ n,k \in {\bf N}$,
denote such $k=l$, when $\psi(x)=\varphi(x)=0$, there exists nontrivial solution
$$
u(x,t)=\left\{
\begin{array}{l}
\left[E_{\alpha,1}\left(-\lambda_lt^\alpha\right)+f_lt^\alpha E_{\alpha,\alpha+1}\left(-\lambda_lt^\alpha\right)\right]  \omega_l(x), t\ge0,\\
\left[\left(\frac{f_l}{\sqrt{\lambda_l}}-\sqrt{\lambda_l}\right)\sin l\pi t +\left(1-\frac{f_l}{\lambda_l}\right)\cos l\pi t +\frac{f_l}{\lambda_l}\right] \omega_l(x), t\le0,\\
\end{array}
\right.
$$
where \\
$$
f_l=-\frac{E_{\alpha,1}\left(-\lambda_lq^\alpha\right)}{q^\alpha E_{\alpha,\alpha+1}\left(-\lambda_lq^\alpha\right)}.
$$
For large $k$ we have $p\approx \frac{2n\pi}{\sqrt{\lambda_k}}$ or $p\approx \frac{(2n+1)\pi}{\sqrt{\lambda_k}}$,$\,n\in {\bf N}$ and when $k\rightarrow \infty$ irregular points $p$ is dense in ${\bf R}^+$.
Generally, since $g$ is any positive bounded function, the irregular points can be set containing both infinite irrational and rational points, but we have the following lemma which shows as $k$ is large, those irregular points may concentrate on irrational points.\\
{\bf Lemma 3.1.} For all sufficient large $k$, then for all $p\in Q^+$, there exists $\delta >0$, such that $|\Delta_k|\ge \delta >0$.

{\it Proof.} We can rewrite (3.22) as follows:
$$\Delta_k=E_{\alpha,1}\left(-\lambda_kq^\alpha\right)-\sqrt{\lambda_k+1} \sin (\sqrt{\lambda_k}p+\gamma_k),\eqno (3.26)
$$
where $\displaystyle{\gamma_k=\arcsin \frac{1}{\sqrt{\lambda_k+1}}}$.\\

For large $k$, as $\lambda_k\to \infty$,
$$
E_{\alpha,1}\left(-\lambda_kq^\alpha\right)=
\frac{C}{\lambda_kq^\alpha \Gamma(1-\alpha)}+O\left(\frac{1}{\lambda_k^2}\right)\to 0,  \ \   q>0.\ \   \ \    \eqno (3.27)
$$
Considering asymptotic behavior of $\lambda_k$, we have
$$
l_k=\sqrt{\lambda_k+1}\sin\left(\sqrt{\lambda_k}p+\gamma_k\right)=\sqrt{\lambda_k+1}\sin\left[pk\pi +\frac{p(c+c_k)}{2k\pi}+\gamma_k\right].
$$
Here $c=\int\limits_0^1g(x)dx$ and $c_k$ is a constant depends on $k,\pi$, which tends to zero, when $k\rightarrow+\infty$.

If $p=\frac{m}{n}\in Q^+,\,(m,n)=1$,

(i) if $km=ln$ for some $l=1,2,...$, then
$$
l_k=\left|\sqrt{\lambda_k+1}\sin\left[pk\pi +\frac{p(c+c_k)}{2k\pi}+\gamma_k\right]\right|\sim \left|\sqrt{\lambda_k+1}\left(\frac{pc}{2k\pi}+\frac{1}{\sqrt{\lambda_k+1}}\right)\right|\sim\frac{pc}{2}+1>0.
$$

(ii) if $km\ne ln$ for any $l=1,2,...$, e.g. $km=ln+s,\,1\le s\le n-1$. Since $\left(\frac{p(c+c_k)}{2\pi k}+\gamma_k\right)\rightarrow 0$, there exists $\varepsilon>0$ such that $\left|\sin\left[pk\pi +\frac{p(c+c_k)}{2k\pi}+\gamma_k\right]\right|\ge \varepsilon>0$. So $|l_k|\rightarrow+\infty$.

Hence, according to (3.26) and (3.27)  we deduce the fact that $\Delta_k$ is bounded below by some positive constant for $p\in Q^+$. \\
The Lemma 3.1. is proved.

\subsection{Existence and stability of the solution}

In the following, assume $g(x)\in L^\infty(0,1)$ is positive, $\varphi(x),\,\psi(x)\in H_0^4(0,1)$. By regularity theorem, $\omega_k(x)\in \mathcal{D}^*(\mathcal{L}^*)$.

The following equation
$$
\omega_k''-g\omega_k=-\lambda_k\omega_k,\,\,\forall k\ge 1,\,\omega_k\in \mathcal{D}^*(\mathcal{L}^*) \eqno (3.28)
$$
is valid in $L^2$ sense.  \\
{\bf(i)} \underline{Proof of $f(x)\in L^2 (0,1)$}.

Introduce the following functions:
$$
f_1(x)=\sum_{k=1}^\infty \lambda_k\varphi_k\omega_k(x), \eqno (3.29)
$$
$$
f_2(x)=\sum_{k=1}^\infty \left[-\frac{\lambda_k(\varphi_k-\psi_k)}{\Delta_k}E_{\alpha,1}\left(-\lambda_kq^\alpha\right)\right]\omega_k(x).\eqno (3.30)
$$
Taking (3.19), (3.28) into account from (3.29) we deduce
$$
f_1(x)=\sum_{k=1}^\infty (\varphi g-\varphi'',\omega_k)\omega_k(x).
$$
Since $g\in L^\infty$, we have $\varphi g\in L^2$ and $\left\|\varphi''\right\|_{L^2}\le \left\|\varphi''\right\|_{H_0^4}<+\infty$. Hence $\left(\varphi g-\varphi''\right)\in L^2$, which yields
$$
\left\|f_1(x)\right\|_{L^2}=\left\|\varphi g-\varphi''\right\|_{L^2}.\eqno (3.31)
$$
Later on we will designate by $C$ any constant, since we are not interested in exact values of them.

Considering Theorem 2.2 and Lemma 3.1, from (3.30) we obtain
$$
\left\|f_2(x)\right\|_{L^2}\le \frac{C}{\delta^2}\left(\left\|\varphi(x)\right\|_{L^2}+\left\|\psi(x)\right\|_{L^2}\right).\eqno (3.32)
$$
Taking (3.25), (3.29)-(3.32) and triangle inequality into account, we deduce
$$
\left\|f(x)\right\|_{L^2}\le \left\|f_1(x)\right\|_{L^2}+\left\|f_2(x)\right\|_{L^2}\le C\left(1+\left\|g(x)\right\|_{L^\infty}\right)\left(\left\|\varphi(x)\right\|_{H_0^4}+\left\|\psi(x)\right\|_{H_0^4}\right).
$$
{\bf(ii)} \underline{Some uniform estimation on the right hand of $u(x,t)$ in (3.23)}.\\
We have by embedding theorem $\left\|\omega_k(x)\right\|_{C[0,1]}\le C\left\|\omega_k(x)\right\|_{H_0^1(0,1)}.$
Thus
$$
\left\|\omega_k(x)\right\|_{C[0,1]}\le C \left\|\omega_k(x)\right\|_{H_0^1(0,1)}\le C\left(\left\|\omega_k'(x)\right\|_{L^2[0,1]}+\left\|\omega_k(x)\right\|_{L^2[0,1]}\right).
$$
According to asymptotic behaviors of $\lambda_k$ and $\omega_k$, we have
$$
\left\|\omega_k(x)\right\|_{C[0,1]}\le C\left(\sqrt{\lambda_k}+1\right).
$$
Introduce
$$
A\equiv \sum_{k=1}^\infty \max_{x\in[0,1]}\left|\frac{\varphi_k-\psi_k}{\Delta_k}\left(E_{\alpha,1}\left(-\lambda_kt^\alpha\right)-E_{\alpha,1}\left(-\lambda_kq^\alpha\right)\right)+\varphi_k\right|C\left(\sqrt{\lambda_k}+1\right).
$$
Easy to deduce that
$$
\varphi_k=\left(\varphi,\omega_k\right)=-\frac{1}{\lambda_k^2}\left(g^2\varphi-2g'\varphi'-2g\varphi''+g''\varphi+\varphi^{IV},\omega_k\right).
$$
Since $g\in L^\infty(0,1)$ and considering
$$
\left\|g^2\varphi\right\|_{L^2}\le \left\|g\right\|_{L^\infty}^2\left\|\varphi\right\|_{H_0^4},\,
\left\|g'\varphi'\right\|_{L^2}\le \left\|g\right\|_{H^{-1}}\left\|\varphi'\right\|_{H^1}\le C\left\|g\right\|_{L^\infty}\left\|\varphi\right\|_{H_0^4},\,
$$
$$
\left\|g\varphi''\right\|_{L^2}\le \left\|g\right\|_{L^\infty}\left\|\varphi\right\|_{H_0^4},\,
\left\|g''\varphi\right\|_{L^2}\le C\left\|g\right\|_{L^\infty}\left\|\varphi\right\|_{H_0^4},\,
\left\|\varphi^{IV}\right\|_{L^2}\le\left\|\varphi\right\|_{H_0^4}
$$
and designating $
G(\varphi)\equiv g^2-2g'\varphi'-2g\varphi''+g''\varphi+\varphi^{IV}
$,
we obtain
$$
\left\|G(\varphi)\right\|_{L^2}\le C\left(1+\left\|g\right\|_{L^\infty}+\left\|g\right\|_{L^\infty}^2\right)\left\|\varphi\right\|_{H_0^4}.
$$
Similarly, bearing in mind $\varphi_k-\psi_k=-\frac{1}{\lambda_k^2}\left(G(\varphi-\psi),\omega_k\right)$, where
$$
G(\varphi-\psi)\equiv g^2(\varphi-\psi)-2g'(\varphi-\psi)'-2g(\varphi-\psi)''+g''(\varphi-\psi)+(\varphi-\psi)^{IV},
$$
we have $G(\varphi-\psi)\in L^2$. Taking all these and Theorem 2.2 into account, we obtain
$$
A\le C\left(\left\|G(\varphi-\psi)\right\|_{L^2}+\left\|G(\varphi)\right\|_{L^2}\right)^{1/2}\sqrt{\sum_{k=1}^\infty\frac{(\sqrt{\lambda_k}+1)^2}{\lambda_k^4}}.
$$
Since $\lambda_k\sim k^2\pi^2+O(1)$, we can state that $A<+\infty$.\\
Introducing
$$
B\equiv \sum_{k=1}^\infty\max_{x\in[0,1]}\left|-\frac{\lambda_k(\varphi_k-\psi_k)}{\Delta_k}t^{\alpha-1}E_{\alpha,\alpha}\left(-\lambda_kt^\alpha\right)\right|C\left(\sqrt{\lambda_k}+1\right)
$$
and by similar evaluations as above, we get
$$
B\le C\left(\left\|G(\varphi-\psi)\right\|_{L^2}\right)^{1/2}\sqrt{\sum_{k=1}^\infty\frac{(\sqrt{\lambda_k}+1)^2}{\lambda_k^2}}
<+\infty.$$
The series in (3.23) converge uniformly in $x\in[0,1]$ and $t\in[0,q]$ and could be differentiated part by part with respect to $t$. One can similarly prove the same fact for $t\in [-p,0]$.\\

{\bf Remark 3.2.} In the above,   $\varphi(x),\,\psi(x)\in H_0^2(0,1)$ are actually enough, if we make a slight different proof, the real problem will happen in a similar proof in (3.24), omitted. There we must need $\varphi(x),\,\psi(x)\in H_0^4(0,1)$, which reflects the bad regularity in the hyperbolic part.\\
{\bf(iii)} \underline{Proof of $u(x,t)\in C\left([-p,q];L^2(0,1)\right)$ and $\displaystyle{\frac{\partial u(\cdot,t)}{\partial t}\in C\left([-p,0]\cup (0,q];L^2(0,1)\right)}$ }.\\
It is easy to verify that
$$
\left\|u(x,t)\right\|_{L^2}\le \frac{C}{\delta^2}\left(\left\|\varphi\right\|_{L^2}+\left\|\psi\right\|_{L^2}\right)+\left\|\varphi\right\|_{L^2},\,\,t\ge 0.
$$
and furthermore, we have $u(x,t)\in C\left([-p,q];L^2(0,1)\right)$.\\ Because of \textbf{(ii)}, $\displaystyle{\frac{\partial u(\cdot,t)}{\partial t}}$ exists and is equal to $U(\cdot,t)\in L^2(0,1)$, where
$$
U(\cdot,t)=\sum_{k=1}^\infty(-\lambda_k)t^{\alpha-1}E_{\alpha,\alpha}\left(-\lambda_kt^\alpha\right)\omega_k(x).
$$ Since, $B<+\infty$ and according to the Theorem 2.2,
$$
\left|(t+h)^{\alpha-1}E_{\alpha,\alpha}\left(-\lambda_k(t+h)^\alpha\right)-t^{\alpha-1}E_{\alpha,\alpha}\left(-\lambda_kt^\alpha\right)\right|   \eqno (3.33)
$$
is bounded for fixed $t>0,\, \forall k\in{\bf N}$. By using the Lebesgue convergence theorem, one can easily verify $\displaystyle{\frac{\partial u(\cdot,t)}{\partial t}\in C\left((0,q];L^2(0,1)\right)}$  Similarly we can prove that\\
$\displaystyle{\frac{\partial u}{\partial t}(\cdot,t)\in C\left([-p,0];L^2(0,1)\right)}.$\\

{\bf Remark 3.3.} In the above proof, it is natural that $\displaystyle{\frac{\partial u}{\partial t}(\cdot,t)}$ is not continuous at $t=0$, since we let $0<\alpha<1$, which is the fractional case, if $\alpha=1$, the equation is just classical parabolic type, and there is no singularity in (3.33), then  $\displaystyle{\frac{\partial u(\cdot,t)}{\partial t}\in C\left([-p,q];L^2(0,1)\right)}$ and the below proof is not needed. In the below, we only consider that $0<\alpha<1$.

{\bf (iv)} \underline{Proof of ${}_CD_{0t}^\alpha u(\cdot,t)\in C\left([0,q];L^2(0,1)\right)$}.\\

Using the definition of the Caputo derivative,
$$
\left\|{}_CD_{0t}^\alpha u(\cdot,t)\right\|_{L^2}= \left\|\int\limits_0^t \frac{s^{\alpha-1}(t-s)^{-\alpha}}{\Gamma(1-\alpha)}\left[\sum_{k=0}^\infty\left\{\frac{\lambda_k(\varphi_k-\psi_k)}{\Delta_k}E_{\alpha,\alpha}\left(-\lambda_ks^\alpha\right)\right\}\omega_k(\cdot)\right]ds\right\|_{L^2}.
$$
According to \textbf{(ii)}, we have
$$
\left\|\sum_{k=1}^\infty\left\{\frac{\lambda_k(\varphi_k-\psi_k)}{\Delta_k}E_{\alpha,\alpha}\left(-\lambda_ks^\alpha\right)\right\}\omega_k\right\|_{L^\infty}<\infty,
$$
which yields $
\left\|{}_CD_{0t}^\alpha u(\cdot,t)\right\|_{L^2}<\infty,
$ since  $\displaystyle{\int\limits_0^t \frac{s^{\alpha-1}(t-s)^{-\alpha}}{\Gamma(1-\alpha)}ds<\infty.}$

Let $t\ge 0$ be fixed, $ t, t+h\in [0,q]$,
$$
\left\|{}_CD_{0t}^\alpha u(\cdot,t+h)-{}_CD_{0t}^\alpha u(\cdot,t)\right\|_{L^2}= \left\|\int\limits_0^t \frac{s^{\alpha-1}(t-s)^{-\alpha}}{\Gamma(1-\alpha)}\left[\sum_{k=0}^\infty\left\{\frac{\lambda_k(\varphi_k-\psi_k)}{\Delta_k}N_k\right\}\omega_k(\cdot)\right]ds\right\|_{L^2},
$$
where $$N_k=\left|\left(1+\frac{h}{s}\right)^{\alpha-1}E_{\alpha,\alpha}\left(-\lambda_k(s+h)^\alpha\right)-E_{\alpha,\alpha}\left(-\lambda_ks^\alpha\right)\right|.
$$
Since $0<\alpha<1$, so for all $t\in [0,q]$, there exists some $C>0$ such that $\left|(1+\frac{h}{s})^{\alpha-1}\right|\le C$ uniformly, together with the properties of Mittag-Leffler function, in terms of asymptotic behavior of $\lambda_k$, we have $N_k\le C'$ ;  With the similar definition of $A$ , using the Lebesgue convergence theorem,  we have
$$
\lim_{h\to0}\left\|{}_CD_{0t}^\alpha u(\cdot,t+h)-{}_CD_{0t}^\alpha u(\cdot,t)\right\|_{L^2}\to0.
$$
{\bf (v)} \underline{Proof of $ u(\cdot,t)\in C\left([0,q]; \mathcal{D}^*(\mathcal{L}^*)\right)$}.

We will consider the following series
$$
 U_m(x,t)=\sum_{k=1}^{m}d_k\omega_k=\sum_{k=1}^{m}\left\{\frac{\varphi_k-\psi_k}{\Delta_k}\left[E_{\alpha,1}\left(-\lambda_kt^\alpha\right)-E_{\alpha,1}\left(-\lambda_kq^\alpha\right)\right]+\varphi_k\right\}\omega_k(x),
$$
$$
F_m(x)=\sum_{k=1}^{m}\left\{\frac{-\lambda_k\left(\varphi_k-\psi_k\right)}{\Delta_k}E_{\alpha,1}\left(-\lambda_kq^\alpha\right)+\lambda_k\varphi_k\right\}\omega_k(x), \ \ \ \ \  t\in [0,q]
$$
considering that $(F_m,\omega_k)=(f,\omega_k)$,
we have
$$
\left({}_CD_{0t}^\alpha U_m,\omega_k\right)+T\left[U_m,\omega_k;t\right]=(f,\omega_k),\,\,(0\le t\le q,\,k=1,2,...m),\eqno (3.34)
$$
where
$$
T\left[U_m,\omega_k;t\right]=\int\limits_0^1 (U_{m,x}\omega_{k,x}+gU_m\omega_k)dx.
$$
Multiplying both side of (3.34) with $d_k$ and sum from $1$ to $m$, we get
$$
\left({}_CD_{0t}^\alpha U_m,U_m\right)+T\left[U_m,U_m;t\right]=(f,U_m),\,\,(0\le t\le q,\,k=1,2,...m).
$$
Taking
$$
\left|\left({}_CD_{0t}^\alpha U_m,U_m\right)\right|\le\left\|{}_CD_{0t}^\alpha U_m\right\|_{L^2}\left\|U_m\right\|_{L^2},\,\left|(f,U_m)\right|\le \frac{1}{2}\left\|f\right\|_{L^2}^2+\frac{1}{2}\left\|U_m\right\|_{L^2}^2,
$$
Garding's inequality [16, p.292], i.e.
$$
T\left[U_m,\omega_k;t\right]\ge \beta\left\|U_m\right\|_{H_0^1(0,1)}^2-\gamma\left\|U_m\right\|_{L^2(0,1)}^2,\,\,\beta>0,\,\gamma\ge 0,
$$
and as well  $\left\|U_m\right\|_{L^2}\le \left\|u\right\|_{L^2}$, into account, we obtain
$$
\beta\left\|U_m\right\|_{H_0^1(0,1)}^2\le C\left\|f\right\|_{L^2}^2+C\left\|U_m\right\|_{L^2}^2+\left\|{}_CD_{0t}^\alpha U_m\right\|_{L^2}\left\|U_m\right\|_{L^2}.
$$
Thus
$$
\left\|U_m\right\|_{H_0^1(0,1)}\le C\left(1+\left\|g\right\|_{L^\infty}\right)\left(\left\|\varphi\right\|_{H_0^4}+\left\|\psi\right\|_{H_0^4}\right).
$$
$U_m$ is uniformly bounded in $H_0^1(0,1)$ and as well
$$
\left\|U_m\right\|_{L^2\left([0,q];H_0^1(0,1)\right)}\le C, \,\,\left\|U_m\right\|_{L^\infty\left([0,q];H_0^1(0,1)\right)}\le C.
$$

There exists a subsequence $\{U_{ml}\}_{l=1}^\infty\subset \{U_{m}\}_{m=1}^\infty$, and $u\in L^2\left([0,q];H_0^1(0,1)\right)$ such that $U_{m,l}\rightarrow u$ weakly.

By standard approximation arguments we see
$$
\left({}_CD_{0t}^\alpha u,v\right)+T\left[u,v;t\right]=(f,v),\,\,\forall v\in H_0^1.
$$
Above $u$ is unique for all $t\in[0,q]$ and by definition it is the same as (3.23).

We rewrite it as
$$
T\left[u,v\right]=(\bar{h},v),
$$
where $\bar{h}=f-{}_CD_{0t}^\alpha u\in L^2(0,1)$ for all $t\in[0,q]$.

From elliptic regularity theorem [17, p.317], we know $u\in H^2(0,1)$ for $0\le t\le q$ and
$$
\left\|u\right\|_{H^2(0,1)}^2\le C\left(\left\|\bar{h}\right\|_{L^2(0,1)}^2+\left\|u\right\|_{L^2(0,1)}^2\right)\le C\left(\left\|f\right\|_{L^2(0,1)}^2+\left\|{}_CD_{0t}^\alpha u\right\|_{L^2(0,1)}^2+\left\|u\right\|_{L^2(0,1)}^2\right).
$$
Since $u''={}_CD_{0t}^\alpha u-q(x)u-f(x)$, and ${}_CD_{0t}^\alpha u, u$ are continuous with respect to $t$ in $L^2$, so $\displaystyle{\lim_{h\to0}\left\|u(t+h,\cdot)-u(t,\cdot)\right\|_{H^2(0,1)}^2\to 0}$. Thus, $u(x,t)\in C\left([0,q];\mathcal{D}^*(\mathcal{L}^*)\right)$.

Furthermore, by taking into consideration hyperbolic part, we get
$$
u(x,t)\in C\left([-p,q];\mathcal{D}^*(\mathcal{L}^*)\right)
$$
and ${}_CD_{0t}^\alpha u(\cdot,t)\in C\left([-p,0);L^2(0,1)\right)$. Some verifications of (3.10)-(3.14) can be found in Appendix.

In case $\varphi(x)=\psi(x)\in H_0^2(0,1)$ we will get the same result.

{\bf Remark 3.4.} If we replace first condition  of (3.8) with $u(0,t)=u_x(1,t)=0$, we have asymptotic behavior of eigenvalues[15, p.140]
$$
\lambda_n=\left(n+1/2\right)^2\pi^2 +\int\limits_0^1g(x)dx-\int\limits_0^1g(x)\cos(2n+1)\pi xdx+O\left(\frac{1}{n}\right),
$$
problem is again ill-posed and the proof is the same.\\

{\bf Remark 3.5.} In the case $\alpha =1$ , they are classical equations with $\frac{\partial u}{\partial t}$ instead of  ${}_CD_{0t}^\alpha u$ and the result is similar.\\

{\bf Remark 3.6.} We note that result of this section generalize the work [5] in particular case ($\alpha=1, r(x)=1,\,e(x)=0$).\\

\section{Case, when hyperbolic part of the mixed equation is purely time-fractional wave equation}
\subsection{Formulation of a problem}

Consider equation
$$
{f}(x)=\left\{
\begin{array}{l}
{}_CD_{0t}^\alpha {u}-\mathcal{L}^*{u},\,t>0,\\
{}_CD_{t0}^\beta {u}-\mathcal{L}^*{u},\,t<0,\\
\end{array}
\right.\eqno (4.1)
$$
together with condition
$$
\left\{
\begin{array}{l}
{u}(0,t)={u}(1,t)=0,\,-p\le t\le q,\hfill\\
{u}(x,-p)={\psi}(x),\,{u}(x,q)={\varphi}(x),\,0\le x\le 1,\hfill\\
\displaystyle{\lim_{t\rightarrow +0} {_CD_{0t}^\alpha {u}(x,t)}={u}_t(x,-0),\,\,0<x<1,}
\\
\end{array}
\right.\eqno (4.2)
$$
where $0<\alpha < 1,\,1<\beta<2$.

First we define a weak solution as follows:

{\bf Definition 4.1.} We call $u(x,t)$ a \textit{weak solution} to (4.1)-(4.2), if $f(x)\in L^2(0,1)$ and $u(\cdot, t)\in\mathcal{D}^*\left(\mathcal{L}^*\right)\equiv H^2(0,1)\cap H_0^1(0,1)$ for $t\in [-p, q]$
and
$$
\begin{array}{l}
\displaystyle{{u}\in C\left([-p, q];\mathcal{D}^*\left(\mathcal{L}^*\right)\right),\,
\frac{\partial {u}}{\partial t}\in C\left([-p,0]\cup (0,q];L^2(0,1)\right),}\hfill\\
\displaystyle{{}_CD_{0t}^\alpha {u}\in C\left([0,q];L^2(0,1)\right),
{}_CD_{t0}^\beta {u} \in C\left((-p,0);L^2(0,1)\right),}\hfill\\
\end{array}
\eqno (4.3)
$$
$$
\lim_{t\rightarrow -p}\left\|{u}(\cdot,t)-{\psi}\right\|_{\mathcal{D}^*\left(\mathcal{L}^*\right)}=0,\eqno (4.4)
$$
$$
\lim_{t\rightarrow q}\left\|{u}(\cdot,t)-{\varphi}\right\|_{\mathcal{D}^*\left(\mathcal{L}^*\right)}=0,\eqno (4.5)
$$
$$
\lim_{|t|\rightarrow 0}\left\| {_C}D_{0t}^\alpha {u}(\cdot,t)-{u}_t(\cdot,t)\right\|_{L^2}=0, \eqno (4.6)
$$
$$
\left({}_CD_{0t}^\alpha {u},\eta\right)+\left(\mathcal{L}^*{u},\eta\right)=\left({f},\eta\right),\,\, t\in [0,q],\,\forall \eta\in \mathcal{D}^*\left(\mathcal{L}^*\right),\eqno (4.7)
$$
$$
\left({}_CD_{t0}^\beta {u},\eta\right)+\left(\mathcal{L}^*{u},\eta\right)=\left({f},\eta\right),\,\, t\in [-p,0),\,\forall \eta\in \mathcal{D}^*\left(\mathcal{L}^*\right). \eqno (4.8)
$$

{\bf Problem 2.} To find a weak solution ${u}$ for (4.1), (4.2) and as well function ${f}(x)\in L^2\left(0,1\right)$ in the domain $\Omega^*$, satisfying (4.3)-(4.8).

By similar algorithm as in the problem 1, i.e. representing solution ${u}(x,t)$ and ${f}(x)$ by
$$
\begin{array}{l}
\displaystyle{{u}(x,t)=\sum_{k=1}^\infty {V}_k(t)\omega_k(x),\,t\ge 0,}\\
\displaystyle{{u}(x,t)=\sum_{k=1}^\infty {W}_k(t)\omega_k(x),\,t\le 0,}\\
\displaystyle{{f}(x)=\sum_{k=1}^\infty {f}_k(t)\omega_k(x),\, 0\le x\le 1,}\\
\end{array}
$$
we have
$$
{V}_k(t)={V}_k(0)E_{\alpha,1}\left(-\lambda_kt^\alpha\right)+{f}_k t^\alpha E_{\alpha,\alpha+1}\left(-\lambda_kt^\alpha\right),
$$
at $t>0$ and
$$
{W}_k(t)={W}_k E_{\beta,1}\left(-\lambda_k(-t)^\beta\right)+t{W}_k'(0)E_{\beta,2}\left(-\lambda_k(-t)^\beta\right)+{f}_k(-t)^\beta E_{\beta,\beta+1}\left(-\lambda_k(-t)^\beta\right).
$$
Further, instead of (3.22) we obtain
$$
\tilde{\Delta}_k=E_{\alpha,1}\left(-\lambda_k q^\alpha\right)-\left[E_{\beta,1}\left(-\lambda_kp^\beta\right)+\lambda_k p E_{\beta,2}\left(-\lambda_kp^\beta\right)\right]. \eqno (4.9)
$$

Formal solution of the problem 2 has a form
$$
{u}(x,t)=\sum_{k=1}^{\infty}\left\{\frac{{\varphi}_k-{\psi}_k}{\tilde{\Delta}_k}\left[E_{\alpha,1}\left(-\lambda_kt^\alpha\right)-E_{\alpha,1}\left(-\lambda_kq^\alpha\right)\right]+{\varphi}_k\right\}\omega_k(x),\,t\in [0,q],
$$
$$
\begin{array}{l}
\displaystyle{{u}(x,t)=\sum_{k=1}^{\infty}\{{\varphi}_k-\frac{{\varphi}_k-{\psi}_k}{\tilde{\Delta}_k}\left[E_{\alpha,1}\left(-\lambda_kq^\alpha\right)-E_{\beta,1}\left(-\lambda_k(-t)^\beta\right)-\right.} \hfill\\
\displaystyle{\left.-\lambda_k t E_{\beta,2}\left(-\lambda_k(-t)^\beta\right)\right]\}\omega_k(x),\,t\in [-p,0),} \hfill\\
\end{array}
$$
and
$$
{f}(x)=\sum_{k=1}^{\infty}\left\{\frac{-\lambda_k\left({\varphi}_k-{\psi}_k\right)}{\tilde{\Delta}_k}E_{\alpha,1}\left(-\lambda_kq^\alpha\right)+\lambda_k{\varphi}_k\right\}\omega_k(x),
$$
where ${\varphi}_k=({\varphi}, \omega_k)$, ${\psi}_k=({\psi}, \omega_k)$ and $\tilde{\Delta}_k$ is defined by (4.9).

{\bf Theorem 4.1} Let $g(x)\in L^\infty$ is positive, assume $0<\alpha<1$, $1<\beta<2$.
\begin{enumerate}

\item For any $\varphi(x),\,\psi(x)\in H_0^4(0,1)$ and fixed $p\in {\bf R}^+$ if $q$ is sufficient large;\\
then the Problem 2 has unique weak solution and the following inequality is valid:
$$
\left\|u\right\|_{C\left([-p,q];H^2\cap H_0^1\right)}+\left\|{}_CD_{0t}^\alpha u\right\|_{C\left([0,q];L^2\right)}+\left\|{}_CD_{0t}^\beta u\right\|_{C\left([-p,0);L^2\right)}+\left\|f\right\|_{L^2}\le
$$
$$
\le C\left(1+\left\|g\right\|_{L^\infty}+\left\|g\right\|_{L^\infty}^2\right)
\left(\left\|\varphi\right\|_{H_0^4}+\left\|\psi\right\|_{H_0^4}\right).
$$
or
\item If $\varphi(x)=\psi(x)\in H_0^2(0,1)$;
\end{enumerate}

$$
\left\|{u}\right\|_{C\left([-p,q];H^2\cap H_0^1\right)}+\left\|{}_CD_{0t}^\alpha {u}\right\|_{C\left([0,q];L^2\right)}+\left\|{}_CD_{t0}^\beta{u}\right\|_{C\left([-p,0);L^2\right)}+\left\|{f}\right\|_{L^2}\le
$$
$$
\le C\left(1+\left\|{g}\right\|_{L^\infty}+\left\|{g}\right\|_{L^\infty}^2\right)
\left(\left\|{\varphi}\right\|_{H_0^4}+\left\|{\psi}\right\|_{H_0^4}\right),
$$

Proof can be done similarly as in section 3. Instead of the Lemma 3.1 we use the following lemma.

  {\bf Lemma 4.1.}
  For any fixed $p\in R^+$, if $q$ is sufficient large, then there exists some $\delta >0$, such that $\left|\tilde{\Delta_k}\right|> \delta>0$.

 {\it Proof.} Fixed $p\in R^+$, noting that $\lambda_k$ has a uniform positive lower bound and the
 asymptotic properties of Mittag-Leffler funtions (Theorem 2.2), when  $q$ is sufficient large, then  $\left|E_{\alpha,1}\left(-\lambda_k q^\alpha\right)\right|$
 is small enough uniformly for any $k$.

  Since $ \left|E_{\beta,1}\left(-\lambda_1p^\beta\right)+\lambda_k p E_{\beta,2}\left(-\lambda_1p^\beta\right)\right|=C_1>0$ and
  for any $k$,
  $$ \left|E_{\beta,1}\left(-\lambda_kp^\beta\right)+\lambda_k p E_{\beta,2}\left(-\lambda_kp^\beta\right)\right| \neq0,$$
  when $k$ is large, $ \left|E_{\beta,1}\left(-\lambda_kp^\beta\right)+\lambda_k p E_{\beta,2}\left(-\lambda_kp^\beta\right)\right|\approx \frac{1}{p^{\beta-1} \Gamma(2-\beta)}>0$, \\
 it easily sees that there always exists some constant $c>0$ depending on $p$   such that
 $ \left|E_{\beta,1}\left(-\lambda_kp^\beta\right)+\lambda_k p E_{\beta,2}\left(-\lambda_kp^\beta\right)\right|\ge c >0$, which proves the statement. \\

{\bf Remark 4.1.} If $g\in H^1$, according to elliptic regularity theorem $\omega_k(x)\in H^3\cap H_0^1$. The embedding theorem deduce that
$$
\left\|\omega_k(x)\right\|_{C^2[0,1]}\le C\left\|\omega_k(x)\right\|_{H^3\cap H_0^1}.
$$
In this case, we actually can discuss the solution in classical case and get $u\in C\left([-p,q];C^2(0,1)\right)$.

{\bf Remark 4.2.} Assume, $\varphi, \psi$ uniquely determine $(u,g_1,f_1)$ and $(v,g_2,f_2)$. Then we have the following relation between $g_i$ and $f_i$ (i=1,2):
$$
f_1(x)-g_1(x)v=f_2(x)-g_2(x)v,\,\,u\equiv v,
$$
which shows, we cannot uniquely determine $f(x)$ and $g(x)$ at the same time. But $u(t,x)$ does not depend on  $g(x)$.
Moreover, if $g_1\equiv 1$, we get $\varphi, \psi$ from  $f_1,\, u$, then if we know a prior $f_2$ , we can recover $g_2$.

\begin{center}
{\bf Conclusion.}
\end{center}

\begin{enumerate}
\item Hyperbolic part of mixed equation has strong influence to the uniqueness and stability. Precisely, in case $\beta=2$ (Problem 1), there are certain conditions to the $p$, but in purely fractional case of $\beta$, i.e. $1<\beta<2$ (Problem 2), we have uniqueness and stability without any restriction to $p$ as in the Lemma 4.1.

\item If we consider instead of $f(x)$ some function in form of $f(x)h(x,t)$ with known $h(x,t)\in C^2$, which satisfies $|h|\ge \delta>0$, then the problem can be studied similarly.

\item For general $n$-dimensional case with $\mathcal{L}^* $ is  symmetric elliptic operator, the following properties hold,    $$
\left\|\omega_k\right\|_{H^\sigma}\le C\lambda_k^{\sigma/2}\left\|\omega_k\right\|_{L^2},\,\,\sigma=0,1,2,\,\,\,\,\,C^{-1}k^{2/n}\le \lambda_k\le Ck^{2/n}$$ for some C, we will have similar results.
\end{enumerate}

\section{Appendix}
Let us verify first
$$
\left({}_CD_{0t}^\alpha u(\cdot,t),\omega_k(x)\right)={}_CD_{0t}^\alpha \left(u(\cdot,t),\omega_k(x)\right),\,\,0\le t\le q.\eqno (A.1)
$$
Since $\frac{\partial u}{\partial t}\in C\left((0,q]; L^2(0,1)\right)$,
introducing
$$
J_{\epsilon_1,\epsilon_2}u(\cdot,t)=\frac{1}{\Gamma(1-\alpha)}\int\limits_{\epsilon_1}^{t-\epsilon_2} (t-s)^{-\alpha} \frac{\partial u(\cdot,s)}{\partial s}ds,
$$
then $J_{\epsilon_1,\epsilon_2}u(\cdot,t)\in L^2(0,1),\,\,\epsilon_1\le t\le q-\epsilon_2$ as proved in (iv).
Further we have
$$
\left(J_{\epsilon_1,\epsilon_2}u(\cdot,t),\omega_k\right)=\frac{1}{\Gamma(1-\alpha)}\int\limits_{\epsilon_1}^{t-\epsilon_2} (t-s)^{-\alpha} \left(\frac{\partial u(\cdot,s)}{\partial s},\omega_k\right)ds.
$$
For $\epsilon_1,\epsilon_2\rightarrow 0$ we deduce
$$
\left({}_CD_{0t}^\alpha u(\cdot,t),\omega_k\right)=\frac{1}{\Gamma(1-\alpha)}\int\limits_0^t (t-s)^{-\alpha} \left(\frac{\partial u(\cdot,s)}{\partial s},\omega_k\right)ds.
$$
Bearing in mind $\left(\frac{\partial u(\cdot,s)}{\partial s},\omega_k\right)=\frac{\partial}{\partial s}\left(u(\cdot,s),\omega_k\right)$, we obtain (A.1).

\textbf{Verification of (3.12):}

Using formula (2.2) from (3.23) one can easily get
$$
{}_CD_{0t}^\alpha u(x,t)=\sum_{k=1}^\infty \left\{-\frac{\lambda_k(\varphi_k-\psi_k)}{\Delta_k}E_{\alpha,1}\left(-\lambda_kt^\alpha\right)\right\}\omega_k(x),\,t\ge 0
$$
and differentiating (3.24) with respect to $t$, we obtain
$$
\frac{\partial u(x,t)}{\partial t}=\sum_{k=1}^\infty \left\{-\frac{\lambda_k(\varphi_k-\psi_k)}{\Delta_k}\left[\cos\sqrt{\lambda_k}t-\sqrt{\lambda_k}\sin\sqrt{\lambda_k}t\right]\right\}\omega_k(x),\,t\le 0.
$$

According to \textbf{(iii)}, we have
$$
\left\| {_C}D_{0t}^\alpha u(\cdot,t)-u_t(\cdot,t)\right\|_{L^2}^2=\sum_{k=1}^\infty \frac{-\lambda_k(\varphi_k-\psi_k)}{\Delta_k}\left|E_{\alpha,1}\left(-\lambda_kt^\alpha\right)-\cos\sqrt{\lambda_k}t-\sqrt{\lambda_k}\sin\sqrt{\lambda_k}t\right|^2.
$$
Taking $|t|\rightarrow 0$, by the Theorem 2.2, the Lebesgue convergence theorem and under certain regularity conditions to the given functions $\varphi,\,\,\psi$, we verify (3.12).

Verifications of (3.10), (3.11), (3.13) and (3.14) can be done by similar arguments.

\section{Acknowledgement}

Authors would like to thank Professor Zhang Bo for his useful suggestions and fruitful discussions. This works was partially supported by the program \textbf{ "TWAS-CAS visiting fellowships 2012" }.


\smallskip

\begin{thebibliography}{31}

\bibitem{1}{\it Kilbas A.A., Repin O.A.} An analog of the Tricomi problem for a mixed type equation with a partial fractional derivative. Fractional Calculus and Applied Analysis, 13(1) (2010), pp.69-84.\\
\bibitem {2}{\it Berdyshev A.S., Cabada A. and Karimov E.T.} On a non-local boundary problem for a parabolic-hyperbolic equation involving a Riemann-Liouville fractional differential operator. Nonlinear Analysis, 75 (2012), pp.3268-3273.\\
\bibitem{3}{\it Berdyshev A.S., Karimov E.T. and Akhtaeva N.} Boundary value problems with integral gluing conditions for fractional-order mixed-type equation. International Journal of differential Equations, (2011), Article ID 268465.\\
\bibitem{4}{\it Kirane M., Malik S.A.} Determination of an unknown source term and the temperature distribution for the linear heat equation involving fractional derivative in time. Applied Mathematics and Computation, 218 (2011) 163-170.\\
\bibitem{5}{\it Sabitov K.B., Safin E.M.} The Inverse Problem for a Mixed-Type Parabolic-Hyperbolic Equation in a Rectangular Domain. Russian Mathematics (Iz. VUZ), 2010, Vol. 54, No. 4, pp. 48-54.\\
\bibitem{6}{\it Ashyralyev A., Ozdemir Y.} On numerical solutions for hyperbolic--parabolic equations with the multipoint nonlocal boundary condition, Journal of the Franklin Institute (Article in Press).http://dx.doi.org/10.1016/j.jfranklin.2012.08.007\\
\bibitem{7}{\it Li G., Zhang D., Jia X. and Yamamoto M.} Simultaneous inversion for the space-dependent diffusion coefficient and the fractional order in the time-fractional diffusion equation. Inverse Problems 29 (2013) 065014 (36pp). 
\bibitem{8} {\it Sakamoto K., Yamamoto M.} Initial value/boundary value problems for fractional diffusion-wave
equations and applications to some inverse problems. J. Math. Anal. Appl. 382 (2011) 426-447.
\bibitem{9}{\it Jin Cheng, Nakagawa J., Yamamoto M. and Yamazaki T.} Uniqueness in an inverse problem for a
one-dimensional fractional diffusion equation. Inverse Problems 25 (2009) 115002 (16pp).
\bibitem{10}{\it Shashkov, A.G.} System-structural analysis of the heat exchange processes and its application. Moscow, 1983.
\bibitem{11}{\it Pskhu A.V.} Uravneniya v chastnykh proizvodnykh drobnogo poryadka, in: Partial Differential Equations of Fractional Order, Nauka, Moscow, 2005, p. 200 (in Russian).
\bibitem{12}{\it Podlubny, Igor}. Fractional differential equations. An introduction to fractional derivatives, fractional differential equations, to methods of their solution and some of their applications. Mathematics in Science and Engineering, 198. Academic Press, Inc., San Diego, CA, 1999. xxiv+340 pp.
\bibitem{13}{\it Kilbas, Anatoly A.; Srivastava, Hari M.; Trujillo, Juan J.} Theory and applications of fractional differential equations. North-Holland Mathematics Studies, 204. Elsevier Science B.V., Amsterdam, 2006. xvi+523 pp.
\bibitem{14}{\it Gorenflo R., Luchko Y.} An operational method for solving fractional differential equations with the Caputo derivatives. Acta Mathematica Vietnamica, 24(2), 1999, pp. 207-233.
\bibitem{15}{\it Andreas Kirsch}. An Introduction to the Mathematical Theory of Inverse Problems. Springer New York Dordrecht Heidelberg London, 2011 (Second edition). xiv+323 pp.
\bibitem{16}{\it Renardy, Michael and Rogers, Robert C.} An introduction to partial differential equations. Texts in Applied Mathematics 13, New York: Springer-Verlag, 2004 (Second edition ed.). xiv+414 pp.
\bibitem{17}{\it Lawrence C.Evans.} Partial Differential Equations. Graduate Studies in Mathematics, Volume 19, AMS, 1998, xx+662 pp.
\end{thebibliography}
\end{document}